\documentclass{ar-1col}

\usepackage[comma]{natbib}
\usepackage{url,hyperref,xspace,amssymb}
\usepackage
{listlbls}
\newcommand{\ie}{{\em i.e.\/}\xspace}
\newcommand{\eg}{{\em e.g.\/}\xspace}

\newcommand{\etc}{{\em etc.\/}\xspace}

\newcommand{\HIDE}[1]{ }
\newcommand{\COMMENT}[1]{ }

\newcommand{\hugin}{{\sc Hugin}\xspace}

\newcommand{\pc}{{\sc pc}\xspace}

\newcommand{\defeq}{:=}

\newcommand{\half}{\frac{1}{2}}

\newcommand{\eqref}[1]{\mbox{(\ref{eq:#1})}}

\newcommand{\secref}[1]{\mbox{\S$\,$\ref{sec:#1}}}
\newcommand{\partref}[1]{\mbox{Part~\ref{part:#1}}}
\newcommand{\figref}[1]{\mbox{Figure~\ref{fig:#1}}}

\newcommand{\tabref}[1]{\mbox{Table~\ref{tab:#1}}}

\newcommand{\itref}[1]{\mbox{\ref{it:#1}}}

\newcommand{\exref}[1]{\mbox{Example~\ref{ex:#1}}}

\newcommand{\Secref}[1]{\mbox{Section~\ref{sec:#1}}}

\newcommand{\E}{{\mbox{E}}}

\newcommand{\var}{{\mbox{var}}}

\newcommand{\cip}{\mbox{$\perp\!\!\!\perp$}}
\newcommand{\indo}[2]{\mbox{$#1 \,\cip\, #2$}}
\newcommand{\ind}[3]{\mbox{$#1 \, \cip\, #2 \mid #3$}}

\newtheorem{expl}{Example}
\newtheorem{definer}{Definition}

\newtheorem{algor}{Algorithm}

\newtheorem{rem*}{Remark}

\newcommand{\halm}{\hspace*{\fill} $\Box$\par}

\newenvironment{ex}{\begin{expl}\rm}{\halm\end{expl}}

\newcommand{\pr}[1]{{p}(#1)}

 \renewcommand{\theenumi}{(\roman{enumi})}
\renewcommand{\pr}{\mbox{\rm pr}} 

\newcommand{\bX}{\mbox{$\mathbf X$}}
\newcommand{\bM}{\mbox{$\mathbf M$}}

\newcommand{\bY}{\mbox{$\mathbf Y$}}

\newcommand{\idle}{\mbox{$\emptyset$}}

\newcommand{\ice}{\mbox{\rm ICE}\xspace}
\newcommand{\ace}{\mbox{\rm ACE}\xspace}
\newcommand{\sce}{\mbox{\rm SCE}\xspace}

\newcommand{\setto}{\leftarrow}
\newcommand{\late}{\mbox{\rm LATE}}
\renewcommand{\pc}{\mbox{\rm PC}}
\renewcommand{\pr}{\Pr}

\jname{Xxxx. Xxx. Xxx. Xxx.}
\jvol{AA}
\jyear{YYYY}
\doi{10.1146/((please add article doi))}

\begin{document}
\markboth{A. P. Dawid and M. Musio}{Effects of Causes and Causes of
  Effects}

\title{Effects of Causes and Causes of Effects}

\author{A. Philip Dawid$^1$ and Monica Musio$^2$
  
  \affil{$^1$Statistical Laboratory, University of Cambridge, UK;
    email: apd@statslab.cam.ac.uk}

  \affil{$^2$Dipartimento di
    Matematica ed Informatica, Universit\`a degli Studi di Cagliari,
    Italy; email: mmusio@unica.it}}


\begin{abstract}
  We describe and contrast two distinct problem areas for statistical
  causality: studying the likely effects of an intervention (``effects
  of causes''), and studying whether there is a causal link between the
  observed exposure and outcome in an individual case (``causes of
  effects'').  For each of these, we introduce and compare various
  formal frameworks that have been proposed for that purpose,
  including the decision-theoretic approach, structural equations,
  structural and stochastic causal models, and potential outcomes.  It
  is argued that counterfactual concepts are unnecessary for studying
  effects of causes, but are needed for analysing causes of effects.
  They are however subject to a degree of arbitrariness, which can be
  reduced, though not in general eliminated, by taking account of
  additional structure in the problem.
\end{abstract}

\maketitle


\part*{INTRODUCTION}
\label{part:intor}

\section{Overview}
\label{sec:over}
The enterprise of ``statistical causality'' has seen much activity in
recent years, both in its foundational and theoretical aspects, and in
applications.  However it remains rare to draw the distinction
(recognised by \citet{mill}) between two different problem areas
within it: assessing (in individual cases, or in general) the likely
effects of applied or considered interventions---the problem of
``effects of causes'', EoC; and assessing, in an individual case,
whether or not an observed outcome was caused by an earlier
intervention or exposure---the problem of ``causes of effects'', CoE.
Where this distinction is made, it is typically assumed that both
problems can be represented and addressed using a common theoretical
framework, such as the structural causal model of \citet{pearl:book}.

The purpose of the current article is to emphasise the important
logical and technical differences between EoC and CoE problems, and to
explore and compare the various ways in which problems of each kind
can be and have been formulated.  In particular it is argued that
different tools are appropriate for the two different purposes.

In \secref{lawsuit} we introduce the variety of concerns to be
addressed, in the context of a specific law suit.  In \partref{eoc} we
introduce and compare a variety of formalisms that have been proposed
to address ``Effects of Causes''.  \Secref{eocintro} briefly
summarises some philosophical and implementational issues.
\Secref{partcomp} introduces, with examples, the problem of inference
in the presence of an instrumental variable, which is then used
throughout \partref{eoc} as a hook on which to hang the general
discussion.  \Secref{condip} describes purely probabilistic aspects.
Then \secref{dt} introduces the decision-theoretic approach to EoC,
\secref{sem} an approach based on linear models, \secref{npsem} a
nonparametric generalisation of that, and \secref{po} the approach
based on potential outcomes.

We turn to address ``Causes of Effects''in \partref{coe}, for problems
similar to those of \secref{lawsuit}.  \Secref{coe_intro} points to
the need for counterfactual inference, which can not however totally
resolve the ambiguities inherent in such problems.  Two ways of
conducting counterfactual modelling are described in \secref{constr},
which can both be subsumed in the potential outcome approach of
\secref{pocoe}.  In \secref{emp} we consider how empirical data can be
used to inform CoE analysis, but can not totally resolve the inherent
ambiguities.  In \secref{anal} we apply this to address the legal CoE
issues of \secref{lawsuit}, showing how the basic ambiguity, expressed
by interval bounds on the ``probability of causation'', can be refined
when we can observe other variables in the problem.  \Secref{further}
indicates just how limited our CoE analyses have been, and what
difficulties might attend further extension.  

Our concluding remarks summarise some of the lessons to learned from
this review of the different approaches to EoC and CoE.

\section{A law suit}
\label{sec:lawsuit}

In 2014 a class action (``multidistrict litigation'', MDL) was brought
in the United States by more than three thousand women who sued the
pharmaecutical company Pfizer, claiming that they developed (type 2)
diabetes as a result of taking its drug Lipitor (Atorvastatin Calcium)
\citep{lipitor}.  The plaintiffs identified two ``bellwether cases''
of women making such a claim for closer attention.

In order to succeed in such a suit, the plaintiffs would have to
demonstrate, in succession, to the Law's satisfaction, two points
\citep{sef/dlf/apd:socmeth}:
\begin{description}
\item[General causation:] {\em Can\/} Lipitor cause diabetes?\\
\item[Specific causation:] In the individual cases, {\em did\/}
  Lipitor cause their diabetes?
\end{description}
The eventual judgment was in favour of Pfizer.  It was judged that
general causation had not been established for doses 10mg, 20mg and
40mg of the drug, but could be considered for the 80mg dose.  And with
regard to the bellwether cases, it was judged that specific causation
could not be established.

The distinction the court made between the two varieties of causal
question, general and specific, is fundamental, and occurs in many
contexts.  It has various descriptions.  Philosophers talk of ``type''
and ``token'' causation.  Legal scholars talk of ``group'' and
``individual'' causation, and have coined the expression G2i (``Group
to individual'') for the task of arguing from one to the other
\citep{g2i}.  In statistical contexts we may talk of inference about
``the effects of causes'' (EoC), and ``the causes of effects'' (CoE),
which are the designations we mostly use here.

\section{Statistical and causal questions}
\label{sec:questions}
Questions about individual cases can usefully be organised in a
fourfold classification.  We exemplify these for the bellwether case
of Juanita, who is 55 years old, has a total cholesterol of 250 mg/dL,
LDL of 175 mg/dL, HDL of 46 mg/dL, triglycerides of 142 mg/dL, weighs
176 lbs, and has a body mass index (BMI) of 26.37.

\begin{description}
\item[Forecasting:] Juanita has started taking a 80mg dose of Lipitor
  daily.  Is she likely to develop diabetes?
\item[Backcasting:] Juanita has developed diabetes.  Did she take
  Lipitor, and if so in what dose and for how long?
\item[Decision:] Juanita is considering whether to take Lipitor, but
  is worried about developing diabetes.  What should she do?
\item[Attribution:] Juanita took Lipitor 80mg daily for 3 years, and
  developed diabetes.  Was that because she took Lipitor?
\end{description}

While {\em Forecasting\/} and {\em Backcasting\/} are fundamentally
purely statistical exercises, {\em Decision\/} and {\em Attribution\/}
can be classified as ``causal'' questions---the former addressing
``Effects of Causes''
(EoC), and the latter, ``Causes of Effects'' (CoE).\\

\noindent{\bf Forecasting.}  This is an apparently straightforward
statistical task, at least conceptually: we gather high quality data
on individuals sufficiently like Juanita, taking the same treatment,
and observe the proportion going on to develop diabetes.  In practice
this simple recipe will be complicated by non-random sampling of
cases, differences in background characteristics, difficulties
associated with long-term follow-up, censoring by death, \etc, \etc\@
Handling such complications has been a prime focus of statistical
research over many decades, and though the issues raised are very far
from trivial, they raise no new issues of principle.  But we would
also need to argue that the proportion, estimated from the data, of
individuals developing diabetes can be identified with Juanita's
``individual risk''.  While this does raise some subtle philosophical
issues \citep{apd:indrisk}, they can largely be ignored for practical
purposes.\\


\noindent{\bf Backcasting.}  This refers to the task of ``predicting''
uncertain past events on the basis of later observations.  In a
statistical context, this is most typically performed by application
of Bayes's theorem.  Suppose we do not know whether or not Juanita
took the Lipitor, but, as above, have estimated the two ``forward''
forecast probabilities, under each scenario.  We would also need to
assign a prior probability to the event that she did, in fact, take
the drug.  Bayes's theorem supplies the machinery for combining these
ingredients to produce the required ``backward'' probability that she
indeed took Lipitor, on the basis of her having developed diabetes.
Although such Bayesian inferences have, from the very beginning, often
been described as estimating the ``probabilities of causes'', use of
the term ``cause'' here is not really appropriate, since even if we
can conclude that Juanita had taken Lipitor, that might not have been
the cause of her diabetes.

Applying Bayes's theorem is not the only way to conduct backcasting.
More straightforwardly, we could simply collect a sample of
individuals sufficiently like Juanita, confine attention to those who
develop diabetes, and use the proportion of these who had taken
Lipitor to estimate the desired probability for Juanita.  Indeed,
there are circumstances where this simple approach may be preferable
to the Bayesian route \citep{apd:diagn}.\\

\noindent{\bf Decision.} Forecasting is of fundamental importance in
decision analysis.  Suppose Juanita has not yet started taking
Lipitor, and is considering whether or not to do so.  One of her
concerns is whether she will develop diabetes.  She should thus
consider, and compare, how probable this event is under two possible
scenarios: that she does, or that she does not, take the drug.  This
would requires two separate forecasting exercises, and correspondingly
data from two different sets of individuals, according as they do or
do not take Lipitor.

But new difficulties now arise in gathering and using such data.  In
particular, the very treatment desired by such an individual might be
related to her overall health status, and thus affect her risk of
developing diabetes---even were she not to receive that desired
treatment.  In such as case it becomes problematic to disentangle the
effects of {\em desire\/} for treatment and of {\em application\/} of
treatment.  This is an example of the problem of ``confounding'',
which requires careful attention in such cases.\\

\noindent{\bf Attribution.} Questions of forecasting, backcasting and
decision, although beset with many practical difficulties, can all, in
principle at least, be answered directly by means of probabilities
attached to unknown events of interest, probabilities that can be
estimated given suitable data.  However, a question of
attribution---such as ``did taking Lipitor cause Juanita's
diabetes''---is not so readily resolved.  For what is it now that is
unknown?  We know that Juanita took Lipitor, and we know that she
developed diabetes.  There is no unknown event about which we require
inference.  Rather, it the relationship between these events that is
uncertain---was it causal, or not?  Even to understand what we might
mean by such a question is problematic.

We shall consider how to formalise such questions, and explore just
what can be concluded from data about them, in \partref{coe} below.









\part{Effects of Causes}
\label{part:eoc}

\section{Introduction}
\label{sec:eocintro}
\subsection{Causality and agency}
\label{sec:agency}

Philosophers have debated causality for millennia, and have propounded
a large variety of conceptions and approaches.  Statisticians, on the
other hand, had traditionally been reluctant to imbue their inferences
with causal meaning.  But in recent years much more attention has been
given to what we can now term ``statistical causality''.  Particularly
influential have been the contributions of \cite{dbr:jep}, who
promoted a formulation based on ``potential outcomes'', and of
\cite{pearl:book}, based on graphical representations.

Implicit in both these approaches is the idea of a cause as an
intervention applied to a system, in line with the ``agency''
interpretation of causality
\citep{reichenbach,price:bjps91,hausman:book,woodward:book,woodward:sep}.
A main task for statistical causality is to make inference about the
effects of such interventions---that is, understanding the ``effects
of causes''---on the basis of data.  When making use of data, it is
important to distinguish between data generated through
experimentation and purely observational data.

\subsection{Experiment}
\label{sec:exp}

In an experiment, interventions are made on experimental units
according to some known protocol, often involving randomisation, and
their responses measured.  To the extent that the experimental units
and interventions can be regarded as representative of future
interventions on new units, it is in principle straightforward to
infer what effects those interventions will have in future.  ``Design
and analysis of experiments'' is a major enterprise within modern
statistics, involving many subtle and technical considerations, but no
special issues of principle arise.

\subsection{Observation}
\label{sec:obs}

Things are not so straightforward when the data available are purely
observational, and the process whereby treatment interventions were
applied to units is not known.  For example, when choosing between two
treatments, a doctor may have given one preferentially to those
patients he considers sicker.  Then a simple comparison of the
outcomes in the two treatment groups will be misleading, since even if
there is no difference between the treatments, a difference in
outcomes may be seen because of the difference in general health of
the two treatment groups.  This is the problem of ``confounding'',
which prevents us from taking the observational data at face value.
In such a case it may or may not be possible to assess, by more
sophisticated means, genuine causal effects, depending on what is
observed and what assumptions can reasonably be made.  If we know or
can reasonably assume how the doctor behaved, and have data on the
patient characteristics that the doctor used, then we can make
meaningful comparisons and extract causal conclusions; but---in the
absence of further structure or assumptions---this will not be the
case if either of the conditions fails.

Much of the modern enterprise of statistical causality is focused on
addressing this issue of extracting causal conclusions from
observational data.  In order to do so, it will invariably be
necessary to make assumptions, generally untestable in practice, about
the relationship between the behaviours of the ``idle'' observational
system, which generates the observed data, and the same system under a
specified intervention---which is what it wanted, but is not directly
observed.  Such assumptions are sometimes made explicit, and so open
to reasoned scrutiny and debate, but sometimes they remain implicit
and hidden, being taken for granted without critical examination.  The
kind of relationships required can typically be expressed, explicitly
or implicitly, as asserting the equality of certain ingredients in
both idle and interventional circumstances.  While such invariance
properties have sometimes been taken as the very definition of
causality \citep{bhlmann2018invariance}, they can be applied without
any such philosophical commitment.  (Our own philosophical standpoint
remains that based on agency.)

The ``do-calculus'' (\citet[\S3.4]{pearl:book}; see also
\citet[\S9.7]{apd:annrev}) applies to problems that can be modelled by
means of a directed acyclic graph, representing both assumed
conditional independence properties of the observational regime, and
assumed relationships between the observational and interventional
regimes.  For such a case it supplies a complete method for
determining whether a causal estimand of interest can be identified
from observational data, and if so how.


\section{Instrumental variable}
\label{sec:partcomp}
Below we shall introduce, compare and contrast some of the different
statistical formalisms that have been used to model effects of causes.
To be concrete, we shall consider, for each formalism, how it might
model an instrumental variable problem \citep{bowden:book}.  This
involves, in addition to the treatment variable $X$ and response
variable $Y$, a further observed variable $Z$ (the instrument), and an
unobserved variable $U$---all defined for individuals in a study or
larger population.  Typically $Z$ is binary, $X$ and $Y$ are binary or
continuous, and $U$ is multivariate.  Note that in this problem it is
not possible, without imposing still further structure, to identify
the causal effect of $X$ on $Y$ from observational data.

We suppose:

\begin{enumerate}
  \renewcommand{\theenumi}{(\alph{enumi})}
\item
  \label{it:1}
  $U$ is a set of pre-existing characteristics of the individual
\item
  \label{it:0}
  $Z$ is associated with $X$, but not with $U$
\item While $X$ could in principle be assigned externally, in the
  study it was not
\item
  \label{it:A} ``Exclusion restriction'': Given $X$, and the
  individual characteristics $U$, the response $Y$ is unaffected by
  $Z$.  (This vague requirement will be clarified further below).
\end{enumerate}


\begin{ex}
  \label{ex:encour}
  {\bf Encouragement trial}\\
  In an encouragement trial \citep{holland:pathanal}, students are
  randomly assigned to receive encouragement to study.  However, a
  student may not or may not respond to the encouragement.  Here $Z$
  is a binary assignment indicator, taking value 1 for encouragement,
  0 for no encouragement; $X$ is the number of hours the student
  actually studies; $Y$ is the student's score in the final test; and
  $U$ comprises individual characteristics of the student, that may
  affect both $X$ and $Y$; because of randomisation, $Z$ is
  independent of $U$.  We are interested in how a student's choice of
  study hours affects their test score.
\end{ex}

\begin{ex}
  \label{ex:partcomp}
  {\bf Incomplete compliance}\\
  In a medical trial, each patient is randomly assigned to take either
  active treatment ($Z=1$) or placebo ($Z=0$).  However, the patient
  may not comply with the assignment, so that the treatment actually
  taken, $X = 1$ or $0$, may differ from $Z$.  Finally we observe
  whether the patient recovers ($Y=1$) or not.  We allow for possible
  dependence of both $X$ and $Y$ on further unobserved patient
  characteristics, $U$.  Again, randomisation ensures that $Z$ is
  independent of $U$.  We are principally interested in the effect of
  taking the treatment on recovery.
\end{ex}

\begin{ex}
  \label{ex:partcomp2}
  {\bf Availability trial}\\
  In a variation of \exref{partcomp}, $Z=1$ means the active treatment
  is made available to the patient, $Z=0$ that it is not; $X=1$ if the
  treatment is taken, $X=0$ if not.  It is supposed that if the
  treatment is unavailable ($Z=0$) it can not be taken ($X=0$) (though
  it need not be taken when it is available).
\end{ex}

\begin{ex}
  \label{ex:mendel}
  {\bf Mendelian randomisation \citep{katan:lancet}}\\
  Low serum cholesterol level ($X=1$) is thought to be a risk factor
  for cancer ($Y=1$).  Both serum cholesterol and cancer may be
  affected by indicators of lifestyle ($U$).  Possession of the E2
  allele ($Z=1$) of the apolipoprotein E (APOE) gene is known to be
  associated with low serum cholesterol level: this relationship need
  not be causal, but may arise because APOE is in linkage
  disequilibrium with the actual causative gene.  Since ``Nature''
  randomises the APOE allele at birth, and its level is not thought to
  affect lifestyle, $U$ should not be associated with $Z$.  We are
  interested in whether intervening to raise serum cholesterol could
  lower the risk of cancer.
\end{ex}
There are a number of questions we could ask (but not necessarily be
able to answer) in such examples.  In \exref{partcomp}, these might
include:

\begin{enumerate}
\item \label{it:zy} ``What is the probability of recovery for a
  patient who is assigned to active treatment (irrespective of the
  treatment actually taken)?''
\item \label{it:xy} ``What is probability of recovery for a patient
  who (irrespective of assigned treatment) in fact took active
  treatment?''
\item \label{it:zx} ``What is the probability that a patient who
  recovered complied with the assignment?''
\item \label{it:fzy} ``What is the effect on recovery of assignment to
  treatment?''
\item \label{it:fxy} ``What is the causal effect of taking the
  treatment on recovery?''
\end{enumerate}
Questions \itref{zy}, \itref{xy} and \itref{zx} inhabit the lowest
rung, ``seeing'', of Pearl's ``ladder of causation''
\citep{pearl:why}, the first two being instances of ``forecasting'',
and the last of ``backcasting''.  Questions \itref{fzy} and
\itref{fxy} are on the second rung of the ladder, ``doing'', being
instances of ``decision''.  In the sequel we shall mainly be
interested in \itref{fxy}.  Since $Z$ has been randomised, we could
argue that the ``intention to treat'' question \itref{fzy} is
essentially the same as \itref{zy}, which can be straightforwardly
addressed from the study data on $(Z,Y)$.  However, \itref{fxy} can
not readily be answered in the same way as \itref{xy} (an ``as
treated'' analysis), since $X$ has not been randomised, and any
observed association between $X$ and $Y$ might be due to their common
dependence on $U$.

\section{Seeing: Conditional independence}
\label{sec:condip}
We first consider how to express, formally, purely probabilistic
properties of the observational joint distribution of $(X,Y,Z,U)$.
This is all that is required to address forecasting and backcasting
questions such as \itref{zy}--\itref{zx}.  However, it will not be
possible to formulate, let alone solve, causal queries such as
\itref{fxy} in this setting: these live on the second rung, ``doing''.

Specifically, \itref{0} implies that $Z$ is independent of $U$: using
the notation of \citet{apd:CIST}, we write this as
\begin{equation}
  \label{eq:indzu}
  \indo Z U.  
\end{equation}
Here \itref{A} is interpreted as asserting the probabilistic
independence of $Y$ from $Z$, conditional on $X$ and $U$:
\begin{equation}
  \label{eq:indyzxu}
  \ind Y Z {(X,U)}.
\end{equation}

\subsection{Graphical representation}
\label{sec:graphrep}

\begin{figure}[htbp]
  \begin{center}
    \resizebox{3.5in}{!}{\includegraphics{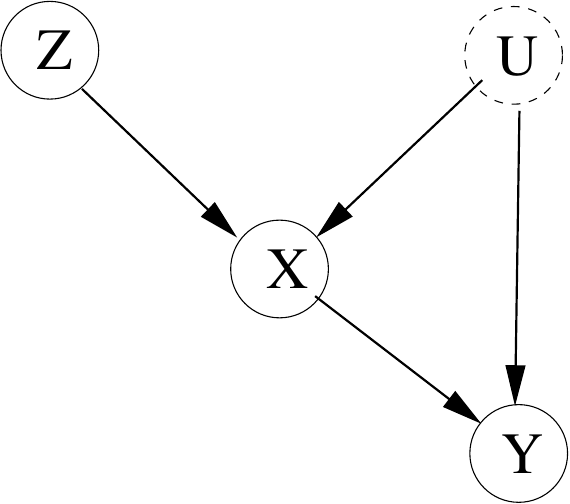}}
    \caption{Instrumental variable: seeing}
    \label{fig:partcomp0}
  \end{center}
\end{figure}

It is often convenient to display such conditional independence
properties by means of a directed acyclic graph (DAG).  Each node in
the DAG represents a variable in the problem, and missing arrows
represent assumed properties of conditional independence in their
joint distribution---see \citet[Section~6]{apd:annrev} for full
details.  The DAG is a partial description, displaying only
qualitative aspects, of the joint distribution.

The DAG representing \eqref{indzu} and \eqref{indyzxu} looks like
\figref{partcomp0} (the dotted outline of $U$ is non-essential, merely
a reminder that $U$ is unobserved).  The absence of an arrow between
$Z$ and $U$ represents their independence, \eqref{indzu}, while the
missing arrow from $Z$ to $Y$ represents their conditional
independence, given the ``parents'' of $Y$, namely $X$ and $U$,
\eqref{indyzxu}.  In general in a DAG representation, any variable is
conditionally independent of its non-descendants, given its parents.
Further conditional independence properties implied by these can be
read off the DAG, using the {\em $d$-separation\/} \citep{verma} or
equivalent {\em moralisation\/} \citep{ldll} criteria, as described in
\citet{apd:annrev}.

The qualitative DAG representation of a joint distribution can be
expanded to a full quantitative description by specifying, for each
variable, its conditional distribution given its parents in the
DAG. (This would be required to encode the condition in \itref{0} that
$X$ is not independent of $Z$).  Elegant algorithms exist, taking
advantage of the DAG structure, for streamlining quantitative
computation of joint and conditional probabilities \citep{mybook}.
Such probabilities are what is needed to address questions of
forecasting and backcasting.

\section{Doing: Decision-theoretic causality}
\label{sec:dt}
The decision-theoretic (DT) approach to causality has been described
in this journal \citep{apd:annrev}; its foundational underpinnings are
examined in \citet{apd:found}.

We have several regimes of interest.  For each possible value $x$ of
$X$ we have an {\em interventional regime\/}, where treatment value
$x$ is forced on an individual (that is, $X$ is ``set'' to $x$, which
we notate as $X \setto x$).  We also have an ``idle'' regime, in which
the treatment $X$ is merely observed, and any value may occur.  It is
helpful to introduce a non-stochastic ``regime indicator'' variable
$F_X$, where $F_X=x$ labels the interventional regime with
$X\setto x$, and $F_X=\emptyset$ labels the idle regime.  The response
variable $Y$ may have different distributions in the different
regimes.  The object of causal inference will usually be some contrast
between the response distributions in the various interventional
regimes---this is what is required to address the decision problem of
choosing which value to set $X$ to.  For example, when $X$ is binary,
interest typically centres on the difference in the expected response
between the two interventional regimes,
$\E(Y \mid F_X=1)-\E(Y \mid F_X=0)$, which is termed the ``average
causal effect'', \ace.

But in the cases to be considered there will be no data available
directly relevant to the interventional settings of interest, and we
shall want to make use of observational data collected under the idle
regime, $F_X=\emptyset$, to make inferences about what would happen in
interventional settings.  This may or may not be possible.  At the
least, it will be necessary to make, and justify, relationships
between the idle and the interventional settings.  DT studies when and
how such relationships can be used to support causal inference from
observational data.

For example, we might be willing to assume, in addition to
\eqref{indzu} and \eqref{indyzxu}, that, {\em no matter whether $X$ is
  merely observed ($F_X=\idle$), or is set by external intervention
  ($F_X=0$ or $1$)}, the following ingredients {\em will be the
  same\/}:
\begin{enumerate}
  \renewcommand{\theenumi}{(\alph{enumi})} \setcounter{enumi}{4}
\item \label{it:zf} the distribution of
  $Z$
\item \label{it:uzf} the distribution of $U$, with $U$ independent of
  $Z$
\item \label{it:yfxu} the conditional distribution of $Y$ given
  $(Z,X,U)$ (which would then in all cases depend only on $(X,U)$,
  since this is so under regime $F_X=\idle$, by \itref{A}.)
\end{enumerate}
These properties do not follow logically from \eqref{indzu} and
\eqref{indyzxu} \citep{apd:wag}, and if they are to be applied they
need additional argument, such as described in \citet{apd:found}.

We can interpret assumptions \itref{zf}--\itref{yfxu} as conditional
independence properties:
\begin{eqnarray}
  \label{eq:indzf}
  Z  &\cip& F_X\\
  \label{eq:indufz}  U  &\cip& (F_X,Z)\\
  \label{eq:indyfxu}
  Y &\cip& (F_X,Z) \mid {(X,U)}.
\end{eqnarray}
Even though $F_X$ is a non-stochastic indicator of regime
(observational/interventional), these intuitively meaningful
``extended conditional independence'' expressions can be manipulated
essentially just as if $F_X$ were a random variable
\citep{pc/apd:eci}.

In this approach, the conditions \itref{zf}--\itref{yfxu}, or
equivalently \eqref{indzf}--\eqref{indyfxu}, which relate to behaviour
under possible intervention at $X$, are the full ``causal
ingredients'' of our model.

\subsection{Graphical representation}
\label{sec:graph}

\begin{figure}[htbp]
  \begin{center}
    \resizebox{1in}{!}{\includegraphics{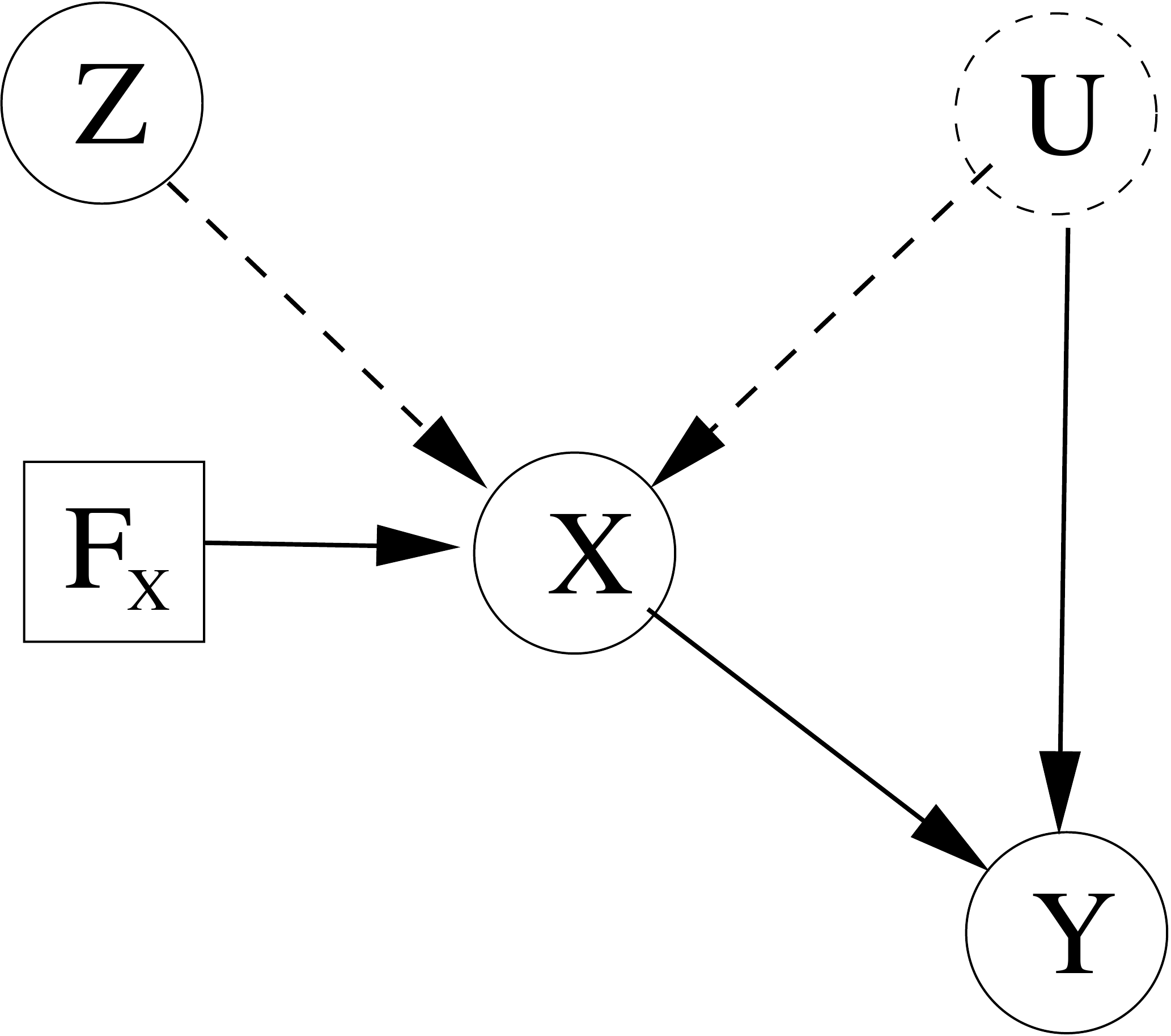}}
    \caption{Instrumental variable: doing}
    \label{fig:partcomp}
  \end{center}
\end{figure}

We can augment \figref{partcomp0} with an additional node for $F_X$
(square to indicate it is non-stochastic).  We obtain
\figref{partcomp}.
This DAG represents (in exactly the same way as before) the assumed
conditional independence assumptions \eqref{indzu}--\eqref{indyfxu},
which fully embody our causal assumptions.  (The dashed arrows from
$Z$ and $U$ to $X$ are there to indicate that they are absent under an
interventional regime $F_X = x$, since then we have $X=x$, trivially
independent of $(Z,U)$).  Note in particular that the arrow from $Z$
to $X$ in \figref{partcomp} does not encode a causal effect of $Z$ on
$X$, since \itref{zf}--\itref{yfxu} are fully consistent with cases,
such as \exref{mendel}, where $Z$ and $X$ are merely associated
\citep[\S10]{apd:beware}.

\subsection{Causal Bayesian Network}
\label{sec:cbn}
\citet{pearl:book} uses the same causal semantics as described above
to construct what he terms a causal Bayesian network (CBN).  The
difference is that he would normally consider the possibility of
intervention on every observable variable---which in our case would
mean adding further intervention indicator nodes $F_Z$, $F_Y$ to
\figref{partcomp}, parents, respectively, of $Z$ and $Y$.  In such a
case the presence of all the intervention nodes is usually taken for
granted, and omitted from the augmented DAG---so rendering it visually
indistinguishable from an unaugmented DAG, here \figref{partcomp0}.
However there are clear advantages to retaining explicit intervention
nodes in the figure:
\begin{enumerate}
\item This eliminates the possibility of confusion between rung~1
  (seeing) and rung~2 (doing) interpretations of apparently identical
  DAGs.
\item The ``causal'' links assumed between regimes are fully
  represented by $d$-separation properties of the augmented DAG.
\item It will (as above) often be appropriate to consider
  interventions on only some of the variables.  In particular, there
  will then be no need to impose the additional cross-regime causal
  constraints associated with further, inessential, intervention
  indicators.
\end{enumerate}

\subsection{Estimation}
\label{sec:est}

Even after assuming links, as above, between the observational and
interventional regimes, it does not follow that we have enough
structure to enable us to use the observational data to estimate, say,
the causal effect, \ace, of $X$ on $Y$.  And indeed, in this example,
further structure must be imposed to support such causal inference.
For instance, in \exref{encour} we might require that $Y$ has a linear
regression on $(X,U)$ (this being the same in all regimes, by
\eqref{indyfxu}):
\begin{equation}
  \label{eq:linreg}
  \E(Y \mid X, U, F_X) = W  + \beta X
\end{equation}
where $W$ is a function of $U$.  Since then
$\E(Y \mid F_X=x) = w_0 + \beta x$, where $w_0 = \E(W)$, $\beta$ has a
clear causal interpretation.  Also, restricting attention to the
observational regime $F_X = \idle$, we obtain
$\E(Y \mid Z) = \E\{\E(Y \mid X, U, Z) \mid Z\} = \E\{\E(Y \mid X, U)
\mid Z\}$, by \eqref{indyfxu}, $= w_0 + \beta\E(X \mid Z)$, by
\eqref{indufz}.  This implies that we can estimate $\beta$, from the
observational data, as the ratio of the coefficients of $Z$ in the
sample linear regressions of $Y$ on $Z$ and of $X$ on $Z$.

In cases such as \exref{partcomp} with binary $X$, \eqref{linreg} is
equivalent to
\begin{equation}
  \label{eq:discreg}
  \sce(u) =  \beta
\end{equation}
a constant for all $u$, where
$\sce(u) = \E(Y \mid U=u, F_X=1)-\E(Y \mid U=u, F_X=0)$ is the {\em
  specific causal effect\/} of $X$ on $Y$, relevant to the
subpopulation having $U=u$.  That is to say, the specific causal
effect is required to be non-random, the same in all subpopulations.
Then $\ace = \E\{\sce(U)\}=\beta$ also, and so is estimable as above.
Alternatively, when all variables are binary, without making any
modelling assumptions we can determine bounds on \ace\ from the data
\citep{b+p:jasa,apd:hsss}.

\section{Linear structural equation model (SEM)}
\label{sec:sem}
Linear structural equation modelling, closely related to path analysis
\citep{wright}, is perhaps the earliest approach to instrumental
variable problems---and much else besides.  It can be considered as an
extension of linear regression modelling.  In the context of the
encouragement trial of \exref{encour}, we might express the
relationship between $Z,X,Y$ by the pair of regression-like equations:
\begin{eqnarray}
  \label{eq:semx}
  X &=& \alpha_0 + \alpha_1 Z + U_X\\
  \label{eq:semy}
  Y &=& \beta_0 + \beta_1 X + U_Y.
\end{eqnarray}
(Such a system would often be completed with a further equation for
$Z$, which here would simply be $Z = U_Z$.  However we omit this on
account of its triviality).  Here $U_X$, $U_Y$ are zero-mean
``residual error'' terms.  In this problem it would be assumed that
$U_X$ and $U_Y$ are uncorrelated with $Z$, but not necessarily with
each other.  The absence of $Z$ in \eqref{semy} embodies the
``exclusion restriction''.

\begin{figure}[htbp]
  \centering \resizebox{2in}{!}{\includegraphics{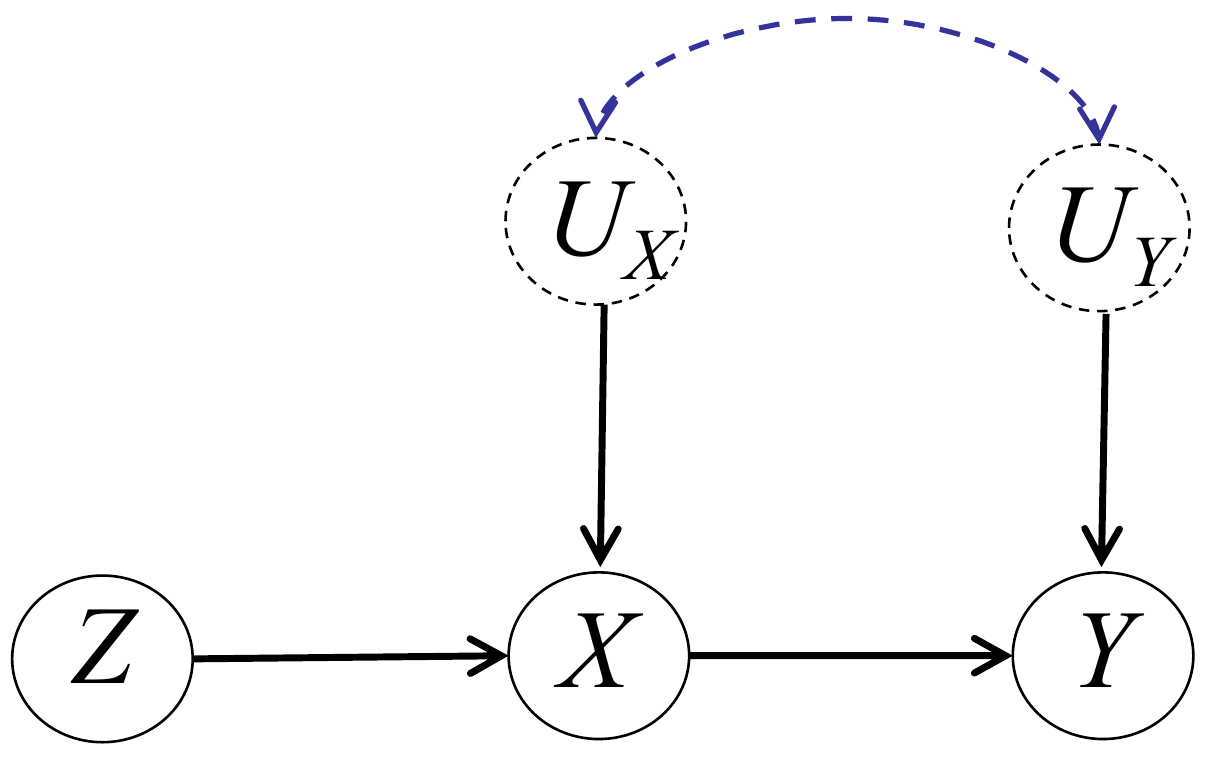}}
  \caption{Structural equation graph}
  \label{fig:sem}
\end{figure}

This model can be rendered graphically as in \figref{sem}.
This may be compared with \figref{partcomp0}, identifying
$U = (U_X, U_Y)$.

As discussed by \citet[\S5.1.2]{pearl:book}, the intended
interpretation---in particular, the causal interpretation---of SEM has
often been unclear.  Pearl's suggestion is as follows.

\begin{enumerate}
\item \label{it:nosetx} In the system \eqref{semx}--\eqref{semy}, $X$
  is functionally determined by $(Z, U_X)$, and $Y$ is functionally
  determined by $(X,U_Y)$.  So we can solve for $(X,Y)$ in terms of
  $(Z, U_X,U_Y)$.  If we have a joint distribution for $(Z, U_X,U_Y)$,
  this determines a joint distribution for $(Z, X, Y)$---which can be
  regarded as representing the undisturbed system.
\item \label{it:setx} If, alternatively, we intervene to set $X$ to
  $x$, it is assumed that we can replace \eqref{semx} by $X=x$, but
  retain \eqref{semy} essentially as is, so that
  $Y = \beta_0 + \beta_1 x + U_Y$, with the distribution of $U_Y$
  unchanged.
\end{enumerate}

This approach gives a causal semantics to a SEM, relating the
observational regime with possible interventional regimes.  As with
any such assumed relationship, it is not to be taken for granted, but
argued for in the context of each particular problem.  When we can
assume \itref{setx}, $\beta_1$ has a clear causal interpretation,
being the rate of change of $\E(Y \mid X \setto x)$ with respect to
the value set, $x$.  However, since $U_Y$ is correlated with $U_X$,
and hence in general with $X$, $\beta$ will not be the coefficient of
$X$ in the observational regression of $Y$ on $X$ in \itref{nosetx},
so can not be identified from that.  Instead we can argue that, since
$U_X$ is uncorrelated with $Z$, $\E(X \mid Z) = \alpha_0 + \alpha_1Z$;
and, since $U_Y$ is uncorrelated with $Z$,
$\E(Y \mid Z) = \beta_0 + \beta_1\E(X \mid Z) = \beta_0 +
\alpha_0\beta_1 + \alpha_1 \beta_1Z$.  It again follows, as in
\secref{est}, that $\beta_1$ can be identified as the ratio of the
coefficients of $Z$ in the observational regressions of $Y$ on $Z$ and
$X$ on $Z$.

\section{Structural Causal Model (SCM)}
\label{sec:npsem}
We can generalise the system \eqref{semx}--\eqref{semy} by dropping
the linearity requirement, yielding
\begin{eqnarray}
  \label{eq:npsemx}
  X &=&  f_X(Z,U_X)\\
  \label{eq:npsemy}
  Y &=& f_Y(X,U_Y)
\end{eqnarray}
where $f_X$, $f_Y$ are specified general functions of their arguments,
and $(U_X, U_Y)$ have a specified joint distribution, typically not
being independent of each other, but being jointly independent of $Z$.
The absence of $Z$ as an argument of $f_Y$ embodies the exclusion
restriction.  \citet{pearl:why} refer to such a nonparametric
structural equation model as a Structural Causal Model (SCM), and we
use this designation in the sequel.

Again, this system can be represented graphically by \figref{sem}.
But with no loss of generality we can use $U$ instead of the pair
$(U_x,U_Y)$, and write the system as
\begin{eqnarray}
  \label{eq:npsemxu}
  X &=&  f_X(Z,U)\\
  \label{eq:npsemyu}
  Y &=& f_Y(X,U)
\end{eqnarray}
with $U$ independent of $Z$.  This system again determines a joint
observational distribution for $(X,Y,Z)$, which is represented by
\figref{partcomp0}.  We might again imbue this structural equation
system with causal semantics (whose relevance in a real life context
will need justification): assume that, under an intervention
$X\setto x$, we can replace \eqref{npsemxu} by $X=x$, and
\eqref{npsemyu} by $Y = f_Y(x,U)$, where $U$ is supposed to retain its
original distribution.  The extended structure is then again
represented by \figref{partcomp}.  In particular, $Y$ will be
independent of $Z$ in an interventional regime, where the dotted
arrows in \figref{partcomp} are absent.  In contrast, in the
observational regime, the distribution of $Y$, given $X=x, Z=z$, is
that of $Y=f_Y(x,U)$ given $f_X(z,U) = x$; because of this
conditioning, the value of $Z$ will typically make a difference to the
distribution of $Y$.

In using \figref{partcomp0} (and, implicitly or explicitly,
\figref{partcomp}) as representations of the SCM system
\eqref{npsemxu}--\eqref{npsemyu}, we are supplying these figures with
yet another semantic interpretation, where the dependence of $X$ and
$Y$ on their parents is taken as deterministic, not stochastic.  This
is to be contrasted with the CBN interpretation of \secref{cbn}, in
which all relationships are allowed to be stochastic.

It can be be shown that, by suitable choice for the distribution of
its $U$ (which distribution is, however, not uniquely determined), the
SCM model can fully reproduce the joint distribution of $(X,Y,Z)$, in
all regimes, implied by a given fully stochastic DT model.  This
property holds in general for any problem represented by a DAG.  For
identifying effects of causes we gain nothing by replacing a
stochastic DT model with a deterministic SCM model.\footnote{See
  however \secref{povars}, and especially \secref{comments}, for
  further discussion of this point.  Also, while \citet{b+p:jasa} make
  essential use of the deterministic functional relationships in the
  SCM to derive estimable bounds on \ace\ in the case of binary
  variables, \citet{apd:hsss} showed how the same bounds can be
  obtained from the purely stochastic DT model.}  In particular, we
again can not identify the causal effect of $X$ on $Y$ without further
assumptions, such as linearity in \eqref{npsemyu}, as for \eqref{semy}
(with $U_Y$ some function of $U$).

\subsection{A comment}
\label{sec:comment}
If, in a SEM or SCM, $U$ is regarded as a persistent attribute of an
individual, the assumed determinism would mean that we would get the
same output each time we applied the same intervention to that
individual.  That would be an unreasonable assumption in most
contexts.  Consequently we should normally consider $U$ as also
incorporating information specific to the occasion of application
(including, perhaps, ``random error''), varying from occasion to
occasion.  Nevertheless, assuming the distribution of $U$ does not
change, average causal effects will still be constant, and so be
meaningful, across occasions.

\section{Potential outcomes}
\label{sec:po}

In the potential outcomes (PO) approach to statistical causality
\citep{dbr:jep}, for each possible value $x$ of the treatment $X$, we
conceive of a version $Y(x)$ of the outcome variable $Y$, all of these
versions co-existing, even before application of treatment.  It is
supposed that $Y(x)$ (or, to be more explicit, $Y(X=x)$) is the
outcome that would be observed in the interventional regime $F_X=x$.
Typically it is further assumed (``consistency'') that in the idle
regime $F_X=\emptyset$ also, whenever $X=x$, the outcome will be
$Y=Y(x)$ (which is why we don't distinguish between $Y(X=x)$ and
$Y(X\setto x)$).  Consistency is required to relate the observational
and interventional regimes.

In the special but common and important case of binary
$X$,\footnote{Of course, similar considerations apply more generally.}
intervening to set $X$ to $1$ would reveal the value of $Y(1)$, while
$Y(0)$ would remain unobserved; and similarly on interchanging 1 and
0.  In this approach, the single response $Y$ is replaced by a
bivariate quantity $\bY = (Y(0),Y(1))$, which must thus be endowed
with a bivariate distribution.  The fundamental causal contrast,
comparing the effects of the two interventions, is considered to be
the {\em individual causal effect\/}, $\ice = Y(1)-Y(0)$.  However,
direct inference about \ice\ is complicated by the fact (termed ``the
fundamental problem of causal inference'' by \cite{pwh:jasa}) that,
because it is logically impossible to intervene on the same individual
in two mutually exclusive ways simultaneously, we can never observe
$\ice$, or estimate its distribution.  For this reason, it is
customary to divert attention to the expected individual causal
effect, $\E(\ice)$.  By linearity of expectation, this is
$\E\{Y(1)\} - \E\{Y(0)\}$, each term of which involves only one
intervention.  This then is the PO version of \ace, as introduced in
\secref{dt}, with essentially the same interpretation.  Note however
that there is no analogue of $\ice$ in the DT approach; neither is
there any DT analogue of, say, $\var(\ice)$, which involves the
correlation between $Y(0)$ and $Y(1)$---a correlation that can never
be estimated, on account of the fundamental problem of causal
inference.

If we start with a SCM representation of a system, we can use it to
construct associated potential outcomes.  For example, starting from
\eqref{npsemyu}, just define, for each $z$, $X(Z=z) = f_X(z,U)$, and
for each $x$, $Y(X=x) = f_Y(x,U)$.  Under the SCM causal semantics,
$X(Z=z) = f_X(z,U)$ is assumed to supply the value of $Y$, when $X=x$,
whether or not there are interventions at $X$ or anywhere else in the
system (except at $Y$ itselvf); this corresponds to the PO consistency
property.  This construction of POs makes them all functions of $U$,
whose distribution thus generates a joint distribution for all POs.

Typically, however, a PO analysis would not make explicit use of an
exogenous variable such as $U$, and might not want to require that
there exist any real-world variable or set of variables $U$ with the
properties assumed in \secref{npsem} (for example, that $Y$ is fully
determined by $X$ and $U$).  Instead, one starts by introducing, as
primitives, jointly distributed stochastic potential outcomes,
$X(Z=z)$, $Y(X=x)$, $Y(Z=z)$, for all possible values of $x$ and $z$,
and work directly with them.  The exclusion restriction now becomes
$Y(Z=z) = Y\{X=X(Z=z)\}$.

Introduce now a new variable $V$, which is simply the collection of
all $X(Z=z)$'s and $Y(X=x)$'s, as $z$ and $x$ vary.  Then $X$ is fully
determined by $(Z,V)$: when $Z=z$, we simply select the relevant
element $X(Z=z)$ of $V$ (this being valid in all regimes, by
consistency); similarly $Y$ is determined by $(X,V)$.  We thus recover
a formal\footnote{In a genuine SCM, $U$ is regarded as a set of
  unobserved real-world background variables, with an appropriate, in
  principle knowable, distribution, that, together with $X$, would
  determine $Y$.  But is is hard to conceive of such a real-world
  interpretation of $V$.} identity with the SCM of \eqref{npsemxu} and
\eqref{npsemyu}, with $V$ substituting for $U$---so long as we have
$V$ independent of $Z$.  If we were starting from a SCM, as above, the
independence of $U$ and $Z$, and thus of $V$ and $Z$, would be easy to
justify, since $U$ represents pre-existing characteristics of the
individual, and $Z$ is randomised.  To make a similar argument for $V$
when taking POs as primitive is more problematic, since $V$ does not
correspond to any real-world quantity (in particular, on account of
``the fundamental problem of causal inference'', certain elements of
$V$, \eg\ $X(Z=1)$ and $X(Z=2)$, are not simultaneously observable).
Nevertheless, the typical assumption is that it is indeed meaningful
to consider the collection $V$ of all possible potential responses as
a pre-existing (albeit unobservable) characteristic of the individual,
and thus argue that $V$ is independent of the randomised variable $Z$.
In this case, we recover a purely formal identity with an SCM model.
Now the linearity condition \eqref{linreg} is equivalent to
$Y(X=x)-Y(X=x') = \beta(x-x')$---which is thus being required to be
non-random.

\subsection{A variation}
\label{sec:povars}
The above specifications can be considered simply as more detailed
ways of realising the CBN structure of \secref{dt}, which is more
general since in a CBN we need not assume the existence of potential
outcomes, which can not be derived from its stochastic form.  And
indeed, for estimating the average causal effect, the extra structure
imposed beyond that of a CBN does not offer any improvement.  But a
SCM or PO approach does allow us to formulate, and purports to solve,
other causal questions.  We consider one such in the context of
\exref{partcomp}, where $Z$, $X$, $Y$, are all binary
\citep{i+a:1994,angrist:96}.

Let ${\bf X}$ denote the pair $(X(Z=0), X(Z=1))$, and $\bf Y$ the pair
$(Y(X=0),(Y(X=1))$.  We shall assume consistency, the exclusion
restriction $Y(Z=z) = Y\{X=X(Z=z)\}$, and
\begin{math}
  \indo Z {({\bf X},{\bf Y})}.
\end{math}
It is easy to see that
\begin{equation}
  \label{eq:yzz}
  Y(Z=z) = Y(X=1)X(Z=z) + Y(X=0)\{1-X(Z=z)\}.
\end{equation}

We can define the following ``individual causal effects'':
\begin{eqnarray}
  \label{eq:icezx}
  \ice_{Z\rightarrow X} &=& X(Z=1)- X(Z=0)\\
  \label{eq:icexz}
  \ice_{Z\rightarrow Y} &=& Y(Z=1)- Y(Z=0)\\
  \label{eq:icexy}
  \ice_{X\rightarrow Y} &=& Y(X=1)- Y(X=0)
\end{eqnarray}
and deduce from \eqref{yzz} that
\begin{equation}
  \label{eq:ices}
  \ice_{Z\rightarrow Y} = \ice_{Z\rightarrow X} \times \ice_{X\rightarrow Y}.
\end{equation}
We note that, since $Z$ is randomised,
$\ace_{Z\rightarrow X} = \E(\ice_{Z\rightarrow X}) = \E(X \mid Z=1) -
\E(X \mid Z=0)$ is readily estimable from the observational data, and
so likewise is
$\ace_{Z\rightarrow Y} = \E(Y \mid Z=1) - \E(Y \mid Z=0)$.  However,
there is no immediate parallel for $\ace_{X\rightarrow Y}$, since $X$
has not been randomised.

If we could replace each \ice\ term in \eqref{ices} by its expectation
\ace, we would have
\begin{equation}
  \label{eq:nogo}
  \ace_{X\rightarrow Y} = \ace_{Z\rightarrow Y} / \ace_{Z\rightarrow X}
\end{equation}
where the right hand-side of \eqref{nogo} is estimable from the
observational data (it is assumed that $Z$ has a causal effect on $X$,
so that $\ace_{Z\rightarrow X}\neq 0$).

When we can assume \eqref{discreg}, $\ice_{X\rightarrow Y}=\beta$ is
non-random, we can take expectations in \eqref{ices}, and \eqref{nogo}
does indeed hold, allowing estimation of $\beta$.  But more generally
${\bf X}$ and ${\bf Y}$ are not independent of each other (when
constructed from a SCM, they both involve the same variable $U$), and
so neither are $\ice_{Z\rightarrow X}$ and $\ice_{X\rightarrow Y}$.
So we can not just take expectations of all terms in \eqref{ices}, and
\eqref{nogo} is typically not valid.

To make further progress, other assumptions must be imposed, in
particular, {\em monotonicity\/}:
\begin{equation}
  \label{eq:mon}
  X(Z=1) \geq X(Z=0).
\end{equation}
That is to say, we do not have any ``defiers'', for which both
$X(Z=0) = 1$ (treatment would be taken when not assigned) and
$X(Z=1) = 0$ (treatment would not be taken when assigned).

Even monotonicity is not sufficient to allow estimation of
$\ace_{X\rightarrow Y}$.  However, it does allow a new interpretation
of the right-hand side of \eqref{nogo}.  For it implies that the
individual causal effect $\ice_{Z\rightarrow X}$ of \eqref{icezx} is
either 1 or 0.  Thus, from \eqref{ices},
\begin{eqnarray*}
  \ace_{Z\rightarrow Y} &=& \E(\ice_{X\rightarrow Y} \mid \ice_{Z\rightarrow X}=1)\times \Pr(\ice_{Z\rightarrow X} = 1)\\
                        &=& \E(\ice_{X\rightarrow Y} \mid \ice_{Z\rightarrow X}=1)\times \E(\ice_{Z\rightarrow X}).
\end{eqnarray*}
It follows that
\begin{equation}
  \label{eq:late}
  \ace_{Z\rightarrow Y} / \ace_{Z\rightarrow X} = \E(\ice_{X\rightarrow Y} \mid \ice_{Z\rightarrow X}=1).
\end{equation}
The right-hand side of \eqref{late} is termed the ``local average
treatment effect'', \late.  Under monotonicity, \late\ is estimable
from the data, since the left-hand side of \eqref{late} is.

\subsubsection{Critical comments}
\label{sec:comments}
\begin{enumerate}
\item \label{it:deter} Considerations similar to \secref{comment}
  suggest that it would typically be appropriate to regard potential
  outcomes, and so individual causal effects, as varying from one
  occasion to another, only their expectations remaining constant.
\item In general the monotonicity assumption is untestable, since
  (under the assumptions of \itref{deter} above) $X(Z=1)$ and $X(Z=0)$
  can not both be observed on the same occasion.  However, it must
  hold in the case of an availability trial as in \exref{partcomp2},
  where necessarily $X(Z=0) = 0$.  Another extreme case where it can
  be inferred is in the presence of a variable $W$, a complete
  mediator between $Z$ and $X$ (so that $X(Z) = X\{W(Z)\}$), where we
  have empirical evidence that, with probability 1, $W(Z=1) = 1$ and
  $X(W=0) = 0$.  If $X(Z=0)=1$, then we deduce $W(Z=0) =1 = W(Z=1)$,
  and so $X(Z=1) = X(Z=0) =1$, and we have no defiers.
\item \late\ is an average causal effect in a subgroup of the
  population: those for whom both $X(Z=0)=0$ and $X(Z=1)=1$.  These
  are termed ``compliers'', since they would take the treatment if
  assigned to do so, and not take it if not assigned (in an
  availability setting, they are those who would take the treatment if
  assigned to do so).  However it is impossible to tell who belongs to
  this subpopulation by knowing only what treatment was assigned and
  what treatment was taken (in an availability trial, an individual
  who was assigned treatment and took it must be a complier, but we
  still can not tell the status of an individual who was not assigned
  treatment).  Indeed, assuming as in \itref{deter} that ${\bf X}$
  will vary from occasion to occasion, so too will the group of
  compliers.  So the relevance of \late\ in practice is debatable.
  Even its definition, relying as it does on a ``cross-world''
  comparison of potential outcomes under both $Z=1$ and $Z=0$, can be
  criticised as essentially metaphysical and unscientific
  \citep{apd/vd:princstrat}.
\item In cases such as \exref{mendel} where $Z$ is not directly causal
  for $X$, the notation $X(Z=z)$ is meaningless, and the above
  analysis can not even get started.
\end{enumerate}

\part{CAUSES OF EFFECTS}
\label{part:coe}

\section{Introduction}
\label{sec:coe_intro}
Let us consider again the initial attribution example: Juanita took
Lipitor 80mg daily for 3 years and developed diabetes. Was that {\em
  because\/} she took Lipitor?

One way of formulating this CoE question is through what the courts
sometimes refer to as the ``but for'' test: is it the case that, {\em
  but for\/} her having taking Lipitor, the diabetes would not have
developed?  This immediately plunges us into counterfactual
considerations.  We know that, in the actual world, the Lipitor was
taken and diabetes developed, and are asked to contrast this with the
outcome in a counterfactual world, in which ({\em counter\/} to the
known {\em facts\/}) the Lipitor had not been taken.  The problem of
course is that the counterfactual world is, by definition,
unobservable, and even its existence---certainly its uniqueness---are
questionable.

Even in deciding on the exact question, choices have to be made.
Juanita took Lipitor 80mg daily for 3 years.  Did she develop diabetes
because she took the 80mg dose (the only one for which the court
accepted general causation), rather than 40mg?  Did it develop because
she took it for 3 years, rather than 2 years?  Each such choice
conjures up a different counterfactual world for comparison with this
one.  We also have choice over what was the observed response: that
she developed diabetes at some point?; that she developed diabetes
within 1 year of stopping Lipitor?  Detailed specification is
obviously important in cases where the response is death: since death
is certain, even in a counterfactual world, we can never say that an
individual would never have died, but for some exposure.\footnote{Even
  in EoC cases, specification of the outcome matters.  In March 2021,
  a few cases of blood clots among individuals who had received the
  Astra-Zeneca vaccine against Covid-19 were observed, and concerns
  were raised about a possible causal connexion, leading to a pause in
  the vaccine's roll-out in some countries.  When it was pointed out
  that the rate of such clots was in fact lower than in the general
  population, attention turned to the few cases of a specific rare
  presentation of the blood clot, cerebral venous sinus thrombosis,
  and whether the vaccine could cause that. (For this outcome too,
  there was no evidence for causation).}

Under the ``but for'' criterion, ``causation'' is understood as the
case that, in the appropriate counterfactual world, where Juanita did
not take the Lipitor (in the same way that she in fact did), she did
not develop diabetes (in the relevant time-frame).  This is
appropriate when the response is all or nothing.  We can also consider
cases with a continuous response, such as time to death, but then it
is not so clear what the focus of our attention should be.  We might
ask, for example, does death occur later, in the relevant
counterfactual world, than it actually did in this world
\citep{sander:ajph}?

Even when our variables have been carefully specified, and the
relevant counterfactual question formulated, it remains unclear just
how to conceive of and structure the counterfactual world of interest.
\citet{dkl:book} develops an approach based on the ``closest possible
world'' to this one, save only for the change to the exposure; but
this only shifts, not solves, the problem.  There appears to be an
unresolvable ambiguity about our counterfactual contrast.

Clearly there are deep philosophical problems, as well as technical
specification issues, besetting any approach to formulating a CoE
problem.  In the sequel we deal only with the case of binary exposure
and outcome variables, denoted by $X$ and $Y$ respectively, assuming
the above specification problems have been addressed.  But there will
still remain some ambiguity about the relevant counterfactual world,
which will be reflected in ambiguity about the answer to the CoE
question.

In \secref{constr} we introduce two approaches to relating the actual
and counterfactual worlds: SCM and CIM.  The former is essentially
deterministic, while the latter allows some stochastic elements.
However both make assumptions that might be regarded as over-strong,
leading to misleadingly precise answers to the CoE question.  In
\secref{pocoe} we show how each of these models can be reformulated in
terms of ``potential outcomes''.  In \secref{emp} we explain how,
taking full account of real-world data on exposure and outcome, this
approach can handle and quantify the remaining ambiguities, by
supplying an appropriate ``interval of ambiguity'' for the probability
of causation, \pc.

We can narrow the interval of ambiguity for an individual case by
deeper understanding of the mechanisms and processes involved
\citep{beyea+sg:99}---even when we can't access the specific details
of these for the individual case at hand.  We develop this theme in
the remaining sections, showing how information about additional
variables can tighten the bounds on \pc.



\section{Counterfactual constructions}
\label{sec:constr}

\subsection{SCM}
\label{sec:npsemcoe}
The approach of \citet{pearl:book} (see also
\citet{pearl:socmeth,sef/dlf/apd:pearlresp}) to CoE is based on SCMs.
In the case of Juanita, this would involve the introduction of an
unobserved exogenous ``background'' variable $U$, and the assumption
that Juanita's diabetes status $Y$ is fully determined by her Lipitor
status $X$ and $U$: $Y=f_Y(X,U)$ (this requires a conception of $U$ as
comprising all other pre-existing quantities that, together with $X$,
would totally determine $Y$---a collection that may not be easy to
comprehend, let alone specify).  In some contexts it might be
appropriate (``ignorability'') to regard $X$ as independent of $U$ in
the observational regime, as would happen, if, for example, $X$ is
generated by a randomising device.  We do not impose this throughout,
and will specify where we do assume it.

The same functional relationship $Y=f_Y(X,U)$ is assumed to hold
(``consistency'') whether or not $X$ is imposed by external
intervention.  To this invariance requirement, familiar from ``effects
of causes'' analysis, we add another, specific to ``causes of
effects'': that the value of the background variable $U$ be the same
in both the factual world, and in the counterfactual world that we
wish to contrast with this one.

To start with, we assume that the function $f_Y$ and the joint
distribution of $(X,U)$ are known.  These unrealistic requirements are
removed in \secref{emp} below.

In the factual world we have observed $X=1$ and $Y=1$, \ie\
$f_Y(1,U)=1$.  We can express the resulting uncertainty about the
value of $U$ by means of its conditional distribution, given
$X=1, f_Y(1,U)=1$ (or, under ignorability, given only $f_Y(1,U)=1$).
We now turn to consider the counterfactual world.  Although $U$ is
supposed to be the same in both worlds (and thus endowed with the
above conditional distribution), $X$ and $Y$ need not be.  We
introduce ``mirror variables'', $X'$ and $Y'$, as their counterfactual
counterparts.  We retain the general structure across worlds, so that
$Y' = f_Y(X',U)$, both in observational and in interventional
counterfactual regimes.

In the counterfactual world, we now consider the effect of an
intervention $X'\setto 0$.  The value of $Y'$ will be $f_Y(0,U)$.
Using the previous conditional distribution of $U$, we obtain the
counterfactual distribution for $Y'$, given the hypothetical
intervention $X'\setto 0$ and the factual knowledge $X=Y=1$.  We can
thus evaluate the ``probability of causation'' as the probability
that, in this distribution, $Y' = 0$.  That is,
\begin{equation}
  \label{eq:pc}
  \pc = \pr(Y'=0 \mid X=1, Y=1, X'\setto 0).
\end{equation}

\subsection{Stochastic Causal Model (StCM)}
\label{sec:cicoe}
A generalisation of the above model, which does not require
deterministic functional relationships, was suggested by
\citet[\S12]{apd:cinfer}, although it has not been developed in
detail.  We again assume that the variable $U$ retains its identity
across the parallel worlds, and introduce the mirror variables $X'$,
$Y'$.  But we now allow the dependence of $Y$ on $(X,U)$ to be given
by a known conditional probability distribution: this allows a more
liberal attitude to the nature of $U$, which can be a perfectly
normal, in principle observable (though typically not observed),
variable.

We assume that the same stochastic relationship also governs the
dependence of $Y'$ on $(X',U)$.  We again require consistency to
relate observational and interventional regimes: the value of $Y$
(resp., $Y'$) would be the same, whether $X$ (resp., $X'$) arose by
intervention or not.  In order to complete the specification, we
regard $Y'$ and $Y$ as conditionally independent, given $(X,X',U)$
(for example, we might consider them as involving ``random noise'',
operating independently across worlds).

Having now a joint distribution for all variables in all worlds, we
can again compute $\pc = \pr(Y'=0 \mid X=1, Y=1, X'\setto 0)$.
Specifically, by conditional independence,
\begin{eqnarray*}
  \pr(Y' =0 \mid X=1, Y=1, X'\setto 0, U) &=& \pr(Y'=0 \mid X'\setto 0, U)\\
                                          &=& \pr(Y'=0 \mid X' = 0, U)\\
                                          &=& \pr(Y=0 \mid X = 0, U).
\end{eqnarray*}
Then $\pc$ is the expectation of this, in the conditional distribution
of $U$ given $X=1, Y=1, X'\setto 0$.  But setting $X'$ does not affect
the joint distribution of $(X,Y,U)$, so we just compute, from Bayes's
theorem, the posterior density
$$p(u \mid X=1,Y=1) \propto  \pr(X=1 \mid u)\,\pr(Y=1 \mid X=1, u)\,p(u),$$
where $p(u)$ is the prior density of $U$ (and the first term on the
right-hand side can be omitted under ignorability).  Finally,
\begin{equation}
  \label{eq:cipc}
  \pc = \frac{\int \pr(Y=0 \mid X=0, u) \,\pr(X=1 \mid u)\,\pr(Y=1 \mid X=1, u)\, p(u)\, du}
  {\int  \pr(X=1 \mid u)\,\pr(Y=1 \mid X=1, u)\, p(u) \,du}.
\end{equation}

\subsection{Twin network}
\label{sec:twin}
\begin{figure}[htbp]
  \centering \resizebox{2in}{!}{\includegraphics{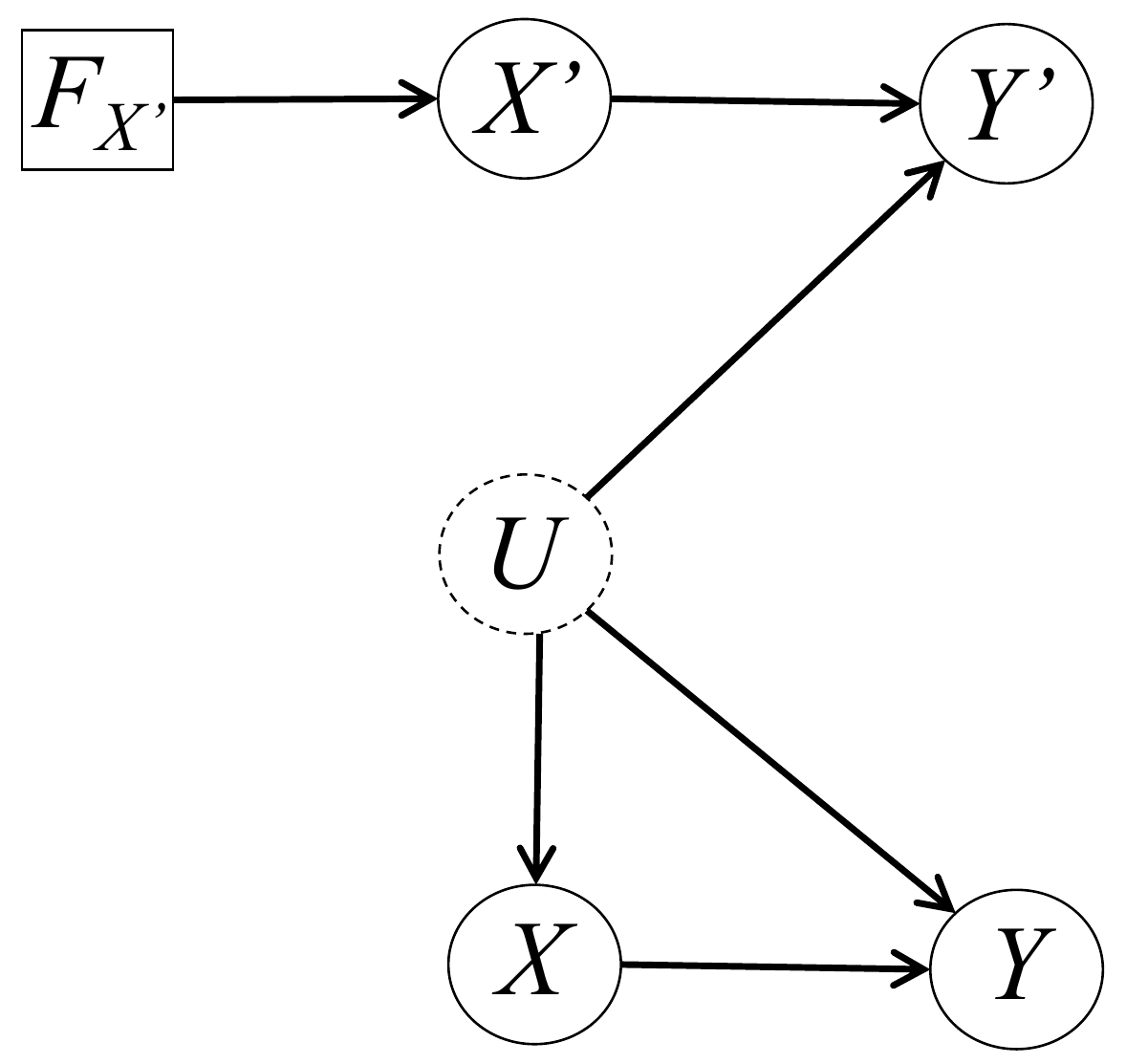}}
  \caption{Twin network.  Under ignorability, the arrow from $U$ to
    $X$ can be removed.  That from $U$ to $X'$ is absent because we
    are considering $X'$ as set by external intervention, taking no
    account of $U$.}
  \label{fig:twin}
\end{figure}

In both the SCM and StCM approaches, the required computation can be
automated by building a ``twin network'' representation of the problem
\citep[\S7.1.4]{pearl:book}, as in \figref{twin}, and making use of
probability propogation algorithms \citep{mybook} as implemented in
software systems such as \hugin.\footnote{\url https://www.hugin.com}
This is useful in more complex problems with more variables.

The factual information $X=1, Y=1$ and the counterfactual intervention
$F_{X'}=0$ (implying $X'=0$) are entered at the relevant nodes, and
propagated through the network to obtain the appropriate conditional
distribution for $Y'$.

\subsection{More general StCMs}
\label{sec:somec}
The StCM approach involves nominating some of the variables in a
problem as shared across worlds, while the others are allowed to
differ.  The associated twin network will have a single copy of the
shared variables and mirror copies of the others, with the original
DAG replicated and stitched together through the shared variables.  As
discussed in \citet{apd:cinfer}, the choice of which variables are to
be regarded as shared\footnote{\citet{mackie:book} refers to these as
  the ``causal field'', and notes the r\^ole of our own choice in its
  specification.} is a matter of imagination rather than science, and
should relate to the specific problem of interest---there can be no
context-free right answer.  For example, there have been law suits by
various states against tobacco companies, claiming that if they had
publicised their knowledge of the dangers of smoking when they first
knew of them, many lives could have been saved.  Damages are sought
for the additional costs placed on health services---meaning the
excess cost in the actual world, over that of an imagined world in
which they had made their knowledge public.  But how should we imagine
that world?  One could reasonably argue that, in such a world, by
giving up smoking, people would have lived longer than they actually
did.  Then the actual \mbox{(non-)}actions of the tobacco companies
might well have {\em increased\/} the cost to the health services.
But what seems to be required for the case at hand is to imagine a
world where people had the same lifetimes, but were healthier, \ie, to
regard lifetimes as shared across parallel worlds---and this even
though lifetimes can be considered an effect of the companies'
decisions.  It is not clear how such considerations could be
accommodated in a SCM.

\section{Potential outcomes}
\label{sec:pocoe}
As observed in \secref{po}, the SCM approach produces implied
potential outcomes, $Y(1) = f_Y(1,U)$, $Y(0) = f_Y(0,U)$.  The pair
$\bY = (Y(1),Y(0))$ is a function of $U$, with a bivariate
distribution induced by that of $U$.  And in fact $\bY$ is all that
needs to be retained of $U$ to fully describe the problem: we can
replace $U$ by $\bY$, with the functional dependence of $Y$ on
$(X,\bY)$ given simply by $Y = Y(X)$.  The problem can thus be more
concisely expressed in terms of $(X, \bY)$, these having a joint
distribution.  Then \eqref{pc} becomes
\begin{equation}
  \label{eq:popc}
  \pc = \pr(Y(0) = 0 \mid X=1, Y(1)=1).
\end{equation}
Note that under ignorability $X$ is independent of $\bY$ (this is
indeed the very definition of ignorability in the PO framework), and
then the conditioning on $X=1$ in \eqref{popc} can be removed.

For the StCM approach we proceed as follows.  The stochastic
dependence of $Y$ on $(X,U)$ can be modelled by introducing a further
unobserved ``noise'' variable $V$, independent of $(X,U)$, and
representing $Y=f_Y(X,U,V)$ for a suitable function $f_Y$.  This can
be done in many ways.  One possible way uses the probability integral
transformation: if $Y$ is a univariate variable whose conditional
distribution function $F_{x,u}(y)$, given $X=x, U=u$, is strictly
increasing, take $V$ to be uniform on $[0,1]$, and
$f_Y(x,u,v) = F^{-1}_{x,u}(v)$.

There will be a counterfactual mirror $V'$ of $V$, with
$Y' = f_Y(X,U,V')$.  We now define potential outcomes
$Y(1) = f_Y(1,U,V)$, $Y(0) = f_Y(0,U,V')$, having a joint distribution
induced by that of $(U,V,V')$.  Although the variables so constructed
will depend on the specific choices made for the noise variable $V$
and the function $f_Y$, it is easy to see that their joint
distribution will in all cases be that of $(Y,Y')$, given
interventions $X\setto 1$, $X'\setto 0$.  And since
$X=1 \Rightarrow Y = Y(1)$, \etc, \eqref{popc} again holds.

We thus see that, in all cases, we can ignore the finer details, and
represent the problem by means of a joint distribution for $(X, \bY)$,
with $\pc$ given by \eqref{popc}.

\section{Empirical information}
\label{sec:emp}
We have so far supposed that the full probabilistic structure of the
model, with its variables $U, X, Y$, is known.  In a StCM, we can take
$U$ to be a specified potentially observable variable, and then this
assumption is not unreasonable.  However it is typically implausible
for a SCM, where, in order to achieve the required deterministic
dependence of $Y$ on $(X,U)$, we would have to conceive of a
fantastically rich $U$.  Alternatively, we can re-express the problem
in terms of the pair $\bY$ of potential responses, with a joint
distribution for $(X,\bY)$.  Without making further assumptions on the
originating SCM or StCM, there are no constraints on this joint
distribution (other than independence of $X$ and $\bY$ under
ignorability).  We can however gather empirical data to constrain it,
and thus hope to estimate $\pc$ by \eqref{popc}.

In the sequel we proceed on this basis, and consider what can indeed
be estimated.  We initially assume that we can only observe $X$ and
$Y$, in interventional and/or observational circumstances.  We can
estimate $\pr\{Y(x) = 1\} = \pr(Y = 1 \mid X \setto x)$ from
interventional studies.  When we can not assume ignorability, we can
also estimate, from observational data,
$\pr\{Y(1) = 1 \mid X = 1\} = \pr(Y=1 \mid X=1)$ (by consistency), and
$\pr\{Y(0) = 1 \mid X = 0\} = \pr(Y=1 \mid X=0)$, as well as the
marginal distribution of $X$.  For this general case it might
initially seem problematic to estimate, say,
$\pr\{Y(1) = 1 \mid X = 0\}$, since this involves non-coexisting
worlds, one with $X=1$ and the other with $X=0$; but we can in fact
solve for it, using
$\pr\{Y(1) = 1\} = \pr\{Y(1) = 1 \mid X = 1\}\,\pr(X=1) + \pr\{Y(1) =1
\mid X = 0\}\,\pr(X=0)$, where all other terms are estimable.

We can thus estimate the bivariate distribution of $(X, Y(1))$, and
likewise that of $(X, Y(0))$.  However, the full trivariate
distribution of $(X, Y(1), Y(0))$ is not estimable: since we can never
observe both $Y(0)$ and $Y(1)$ simultaneously, no data can tell us
directly about the dependence between $Y(0)$ and $Y(1)$, either
marginally or conditionally on $X$.  And since \eqref{popc} requires
such information, typically $\pc$ is not identifiable from empirical
data.

Nevertheless, the estimable bivariate distributions do impose
constraints on the possible values of $\pc$.  Moreover, these
constraints can often be tightened still further when we can observe
other, related, variables in the problem.  We now turn to investigate
such constraints, in a number of contexts.

\section{Analysis of the probability of causation}
\label{sec:anal}
Let us consider the initial attribution example: Juanita took Lipitor
80mg daily for 3 years ($X=1$) and developed diabetes ($Y=1$).  Was
that {\em because\/} she took Lipitor?  We wish to address this
question and assess the probability of causation, \pc, for Juanita's
case, using data collected on other individuals.  To this end we assume
\begin{description}
\item [Exchangeability] Juanita is similar to the population from
  which probabilities have been computed, so that those probabilities
  apply to her.
\end{description}
Exchangeability may require restriction of the data considered to
individuals deemed sufficiently like Juanita.

Except where relaxed in \secref{tian} below, we shall also assume
\begin{description}
\item [Ignorability] The fact that Juanita chose to take the drug is
  not informative about her response to it, either factually or
  counterfactually.  Formally, we require independence,
  $\indo X {\bY}$, between $X$ and the pair of potential responses
  $\bY = (Y(1),Y(0))$.
\end{description} 
Ignorability is a strong assumption, and will often be inappropriate.
When it can be assumed, we can use data from either experimental or
observational studies; otherwise we need data from both of these.

Under ignorability, the target \eqref{popc} becomes
 \begin{equation}
  \label{eq:pcj}
  \pc = \Pr(Y(0) = 0 \mid  Y(1)=1).
\end{equation}
We proceed to assess this using the general potential outcome
framework of \secref{pocoe} and \secref{emp}, where no assumptions are
imposed on the joint distribution of $Y(0)$ and $Y(1)$ beyond those
that can be informed by the empirical data.  Further details may be
found in \citet{dmm:lpr,apd/mm:sef}.

\subsection{Basic inequalities}
\label{sec:base}

Suppose we have access to (observational or experimental) data,
supplying values for
\begin{equation}
  \label{eq:marg}  
  \Pr\{Y(x) = y\} = \Pr(Y=y \mid X= x) \quad (x=0,1
  \quad y=0,1).
\end{equation}
Define
\begin{eqnarray*}
  \tau &\defeq& \Pr(Y=1\mid X= 1) - \Pr(Y=1\mid X= 0)\\
  \rho &\defeq& \Pr(Y=1\mid X= 1) - \Pr(Y=0\mid X= 0). 
\end{eqnarray*}
The joint distribution of $(Y(0),Y(1))$ must have the form of
\tabref{y0y1}, where the marginal probabilities are given by
\eqref{marg}, re-expressed in terms of $\tau$ and $\rho$, and where
the unidentified ``slack'' quantity $\xi$ embodies the residual
ambiguity in the distribution.
\begin{table}[tbp] \centering
  \begin{tabular}[c]{c|cc|c} \multicolumn{1}{c}{}& $Y(1)=0$ &
    \multicolumn{1}{c}{$Y(1)=1$}\\\hline
    $Y(0) = 0$ & $\half(1-\rho - \xi)$ & $\half(\xi+\tau)$ & $\half(1+\tau-\rho)$\\
    $Y(0) = 1$ & $\half (\xi-\tau)$ & $\half(1+\rho-\xi)$ &
                                                            $\half(1-\tau+\rho)$
    \\\hline & $\half(1-\tau-\rho)$ & $\half(1+\tau+\rho)$ & 1
  \end{tabular}
  \caption{Joint distribution of $Y(0)$ and $Y(1)$}
  \label{tab:y0y1}
\end{table}
For all the entries of \tabref{y0y1} to be non-negative, we require
\begin{equation}
  \label{eq:basicineq}
  |\tau| \leq \xi \leq  1-|\rho|.
\end{equation}

The probability of causation  \eqref{pcj} is  
\begin{equation}
  \label{eq:pcP}
  \pc = \frac{\xi+\tau}{1+\tau+\rho}.
\end{equation}
On using \eqref{basicineq}, we obtain the following interval bounds
for $ \mbox{PC} $:
\begin{equation}
  \label{eq:simple}
  l := \max\left\{0, \frac{2\tau}{1+\tau+\rho}\right\} \leq \mbox{PC}  \leq \min\left\{1,\frac{1+\tau-\rho}{1+\tau+\rho}\right\} = :u,
\end{equation}
or equivalently
\begin{equation}
  \label{eq:simpleprobs}
  l = \max\left\{0, 1 - \frac {1}{ \mbox{RR} }  \right\} \leq  \mbox{PC} \leq 
  \min\left\{1,\frac{\Pr(Y = 0 \mid X= 0)}{\Pr(Y = 1 \mid X= 1)}\right\} = u,
\end{equation}
where
\begin{equation}
  \label{eq:rr}
  \mbox{RR}  = \frac{\Pr(Y=1 \mid X= 1)}{\Pr(Y=1 \mid X= 0)}
\end{equation}
is the {\em risk ratio\/} \citep{jr+sg:89}.

In the absence of additional information or assumptions, these bounds
constitute the best available inference regarding $\mbox{PC}$.  In
particular, $\mbox{RR} > 2$, ``doubling the risk'', implies that
$\mbox{PC} > 0.5$.  In a civil legal case, causality might then be
concluded ``on the balance of probabilities''.  However, because of
the remaining ambiguity, expressed by the inequalities in
\eqref{simpleprobs}, finding that $\mbox{RR}$ falls short of 2 does
not imply that $\mbox{PC} < 0.5$.

\subsection{Refining the inequalities: Covariates}
\label{sec:refine}
When we have additional information we may be able to refine our
inferences about $\mbox{PC}$ \citep{kuroki}.

Thus suppose that we also have information on a covariate $S$, a
pretreatment individual characteristic that can vary from person to
person and can have an effect on both $X$ and $Y$.  The relevant
potential responses are now $\bX := (X(s): s\in S)$,
$\bY := (Y(s,x): s\in S, x=0\mbox{ or }1)$ and the relationship
between potential and actual responses is $X = X(S)$, $Y = Y(S,X)$.
We again assume exchangeability and ignorability, the latter now being
formalised as mutual independence between $S$, $\bX$ and $\bY$.

For simplicity we suppose that $S$ is discrete and that we can
estimate from the data the full joint distribution of $(S,X,Y)$.

In the case that we are also able to measure $S$ for Juanita, say
$S = s$, we can simply restrict the experimental subjects to those
having the same covariate value (who are thus like Juanita in all
relevant respects).  The probability of causation is now
$$\pc(s) = \Pr\left(Y(s,0) = 0 \mid Y(s,1) = 1\right).$$
We can bound this just as in \eqref{simpleprobs}, but with all
probabilities now conditioned on $S = s$, obtaining
$l(s) \leq \pc(s)\leq u(s)$.

More interesting is the case in which we don't observe $S$ for
Juanita.  We have to consider what would have been the response if,
counterfactually, Juanita's exposure had been $X=0$.  We assume that
this is the minimal change made between the factual and the
counterfactual worlds, so that, in particular, there is no change to
the value or distribution of $S$.

The probability of causation is now:
\begin{eqnarray}
  \nonumber     
  \mbox{PC} &=& \Pr \left\{Y(S,0) = 0 \mid X(S) = 1, Y(S,1) = 1\right\}\\
            &=& \sum_{s} \mbox{PC}(s) \times \Pr(S= s \mid X=1, Y=1).
\end{eqnarray}

There are no logical relationships between the distributions of
$(Y(s,0),Y(s,1))$ for different values of $S$.  So by independently
varying the values taken by the slack variables in the joint
distribution of these potential responses, all the lower bounds $l(s)$
for $\mbox{PC}(s)$ can be achieved simultaneously.  This leads to an
achievable lower bound for $\mbox{PC}$:
\begin{equation}
  \label{eq:pcl}
  \mbox{PC} \geq L= \sum_{s} l(s) \times \Pr (S= s \mid X=1, Y=1).
\end{equation}
We can express
\begin{eqnarray}
  \nonumber
  L &=& \frac 1 {\Pr(Y=1 \mid X=1)}\\
    &&{}\times
       \sum_{s} \max\left\{0,\Pr(Y=1 \mid X= 1, S= s) - \Pr(Y=1 \mid X= 0, S= s)\right\}\times \Pr(S=s \mid X=1).  \nonumber\\
  \label{eq:nlb}
  \quad
\end{eqnarray}

Similarly we obtain upper bound
\begin{eqnarray}
  \nonumber
  U &=& 1-\frac 1 {\Pr(R=1 \mid E=1)}\\
    &&{}\times
       \sum_{s} \max\{0, \Pr(Y=1 \mid X= 1, S= s) - \Pr(Y=0 \mid X= 0, S= s)\}\times \Pr(S=s \mid X=1).  \nonumber\\
  \label{eq:nub}
  \quad
\end{eqnarray}

We can not compare these bounds directly with those of
\eqref{simpleprobs}, since when we don't take account of $S$ the
relation between $X$ and $Y$ is generally non-ignorable:
$\pr(Y=y \mid X\setto x) \neq \pr(Y=y \mid X= x)$.

\subsection{Non-ignorability}
\label{sec:tian}
\citet{tian/pearl:probcaus} analysed this non-ignorable case where we
don't observe $S$, either for Juanita or in the external data (still
considered exchangeable with Juanita).  We now need both observational
and experimental data on $X$ and $Y$.  \citet{tian/pearl:probcaus}
develop the following lower bound for
$\mbox{PC} = \pr\{Y(0) = 0 \mid X=1, Y(1)=1\}$:
\begin{equation}
  \label{eq:tplb}
  L'= \max\left\{0, \frac{\Pr(Y=1) - \Pr(Y=1 \mid X \leftarrow 0)}{\Pr(X=1, Y=1)}\right\}.
\end{equation}

\citet{apd/mm:sef} show that this can also be derived as a special
case of our expression \eqref{nlb}, if we substitute for $S$ the
binary variable $D=$ ``desired exposure'' \citep{fc/mm:coe}.  $D$ will
be identical with $X$ in an observational context, but need not be so
in an experimental setting, where $D$ may not be observable.

In our case, with access to information on $S$, we could compute
$\Pr(Y=1 \mid X \leftarrow 0)$ by the ``back-door formula'':
\begin{equation}
  \label{eq:backdoor}
  \Pr(Y=1 \mid X \leftarrow 0) = \sum_{s}\Pr(Y=1 \mid X= 0, S= s)\times\Pr(S=s),
\end{equation}
and thus compute $L'$ of \eqref{tplb}.  It can be shown
\citep{apd/mm:sef} that $L' \leq L$ of \eqref{nlb}, with equality if
and only if all the conditional risk ratios
$$\frac{\Pr(Y=1 \mid X = 1, S= s)}{\Pr(Y=1 \mid X = 0, S= s)}\quad\quad (s\in S)$$
lie on the same side of 1: knowing, and using, the information about
$S$ is at least as good as ignoring it.  Similarly we can show that
the upper bound $U$ of \eqref{nub} does not exceed the upper bound
$U'$ derived by \citet{tian/pearl:probcaus},
$$U' = \min\left\{1, \frac{\Pr(Y =0 \mid X \leftarrow 0)-\Pr(X=0,Y=0)}{\Pr(X=1,Y=1)}\right\},$$ 
with equality if and only if all the ratios
$$\frac{\Pr(Y=1 \mid X = 1,S = s)}{\Pr(Y=0 \mid X = 0,S = s)}\quad\quad (s\in S)$$ lie on the same side of 1.

\subsection{Mediators}
\label{sec:mediators}
We now consider the case that a third variable $M$ acts a complete
mediator in the causal pathway $X \rightarrow M \rightarrow Y$ between
the exposure $X$ and the response $Y$.  Again we restrict to the case
that all variables are binary.  We introduce the potential value
$M(x)$ of $M$ for $X= x$, and $Y(m)$, the potential value of $Y$ for
$M= m$, and define $\bM = (M(0), M(1))$, $\bY = (Y(0),Y(1))$.  We
observe $X$, $M=M(X)$ and $Y=Y(M)$.  We assume the exchangeability and
the ignorability conditions, the latter expressed as mutual
independence between $X$, $\bM$ and $\bY$.  This implies the
observational conditional independence
\begin{equation}
  \label{eq:as20}
  \ind {Y} {X} {M}
\end{equation}
which is a testable implication of our assumptions.  We assume we have
data supplying values for $\Pr(M=m \mid X= x)$ and
$\Pr(Y= y \mid M= m)$, and compute, by \eqref{as20},
\begin{equation}
  \label{eq:ab0}
  \Pr(Y=y \mid X = x) = \sum_m\Pr(Y=y \mid M = m)\Pr(M=m \mid
  X = x),
\end{equation}

For the case that $M$ is observed in the experimental data, but not
for Juanita, \citet{apd/rm/mm:sis} showed that this additional
information does not change the lower bound $l$ on $\pc$ in
\eqref{simpleprobs}, but does lower the upper bound, $u$.
\citet{apd/mm:sef} extend this analysis to cases with additional
covariates, while
\citet{dhm} deal with the case that we have a complete mediation
sequence
$X=M_0\rightarrow M_1\rightarrow \cdots\rightarrow M_{n-1}\rightarrow
M_n = Y$, and know the probabilistic structure of each link in the
chain; with all, some, or none of the $M$s being observed for Juanita.

\section{Further CoE problems}
\label{sec:further}
We have only dealt here with the case of a single putative cause,
understanding causation in terms of the ``but for'' criterion.  There
are many more complex problems that can not be handled in this way: in
particular, the whole field of legal causation has to handle a wide
variety of problems involving multiple competing causes and other
concepts of causality \citep{hart:book,goldberg:book}.  While there
have been some interesting statistical treatments of specific
problems, \eg, \citet{lacox}, it seems fair to say that general
philosophical understandings of causality in such problems have not
reached maturity.  To the extent that problems are modelled in formal
terms, this often involves a purely deterministic understanding of
causality, which is not easily translated into a stochastic framework.
\citet{halpern:book} makes an interesting attempt to pin down the
concept of ``the actual cause'' using the SCM framework, but admits
that he is unable to reach a fully satisfying conclusion.

There is clearly much ground remaining to be covered in understanding
``causes of effects'', but it is perhaps premature to attempt more
detailed statistical treatment before clearer general principles have
emerged.

\part*{CONCLUSION}
\label{part:conc}
We have presented a thorough account of a number of ways in which the
statistical problems of effects of causes (EoC) and of causes of
effects (CoE) have been formulated.  Although most treatments of
statistical causality use essentially identical tools to address both
these problem areas, we consider that this is inappropriate.  Popular
formalisms, such as potential outcomes and structural causal models,
involve deterministic relations, and allow formal statements
concerning two or more parallel worlds simultaneously.  While
something of this nature appears unavoidable for CoE considerations,
it is unnecessary for EoC analyses, which can proceed using stochastic
models and statistical decision theory.  Furthermore, the use of an
inappropriate formal framework brings with it the danger that we treat
any mathematically well-formed formula (such as those describing
``individual causal effect'' and ``local average treatment effect'')
as meaningful, when it may not be.

\section*{DISCLOSURE STATEMENT}
The authors are not aware of any affiliations, memberships, funding,
or financial holdings that might be perceived as affecting the
objectivity of this review.

\section*{ACKNOWLEDGMENTS}
The second author was partially supported by the project STAGE of
Fondazione di Sardegna.

\bibliographystyle{ar-style1}
\bibliography{strings,causal}


\end{document}